\documentclass{amsart}

\usepackage{xypic}
 \input xy
\xyoption{all}
\usepackage{epsfig}
\usepackage{amsthm}
\usepackage{amssymb}
\usepackage{amsmath}
\usepackage{amscd}
\usepackage{color}

\newcommand{\bg}{\begin{equation}}
\newcommand{\ed}{\end{equation}}
\newcommand{\bga}{\begin{eqnarray}}
\newcommand{\eda}{\end{eqnarray}}
\newcommand{\pf}{\textbf{Proof:\ }}

\def\cbdu{\par{\raggedleft$\Box$\par}}

\newtheorem {Theorem}  {Theorem}

\numberwithin{Theorem}{section}

\newtheorem {Lemma}[Theorem]  {Lemma}

\theoremstyle{definition}

\theoremstyle{remark}

\expandafter\chardef\csname pre amssym.def
at\endcsname=\the\catcode`\@ \catcode`\@=11
\def\undefine#1{\let#1\undefined}
\def\newsymbol#1#2#3#4#5{\let\next@\relax
 \ifnum#2=\@ne\let\next@\msafam@\else
 \ifnum#2=\tw@\let\next@\msbfam@\fi\fi
 \mathchardef#1="#3\next@#4#5}
\def\mathhexbox@#1#2#3{\relax
 \ifmmode\mathpalette{}{\m@th\mathchar"#1#2#3}%
 \else\leavevmode\hbox{$\m@th\mathchar"#1#2#3$}\fi}
\def\hexnumber@#1{\ifcase#1 0\or 1\or 2\or 3\or 4\or 5\or 6\or 7\or 8\or
 9\or A\or B\or C\or D\or E\or F\fi}

\font\teneufm=eufm10 \font\seveneufm=eufm7 \font\fiveeufm=eufm5
\newfam\eufmfam
\textfont\eufmfam=\teneufm \scriptfont\eufmfam=\seveneufm
\scriptscriptfont\eufmfam=\fiveeufm

\catcode`\@=\csname pre amssym.def at\endcsname

\newcounter{remark}
\usepackage{color}

\newcommand{\R}{\mathbf{R}}

\def  \R   {{\mathbb R}}

\def  \12  {{\frac{1}{2}}}

\def\build#1_#2^#3{\mathrel{\mathop{\kern 0pt#1}\limits_{#2}^{#3}}}

\begin{document}

\title[Norm-inflation for BBM Equation]{Norm-inflation Results for the BBM Equation}

\author[Jerry Bona]{ Jerry Bona}
\address{Department of Mathematics, University of Illinois at Chicago, Chicago, IL 60607,USA}
\email{bona@math.uic.edu}

\author[Mimi Dai]{ Mimi Dai}
\address{Department of Mathematics, University of Illinois at Chicago, Chicago, IL 60607,}
\email{mdai@uic.edu}

\thanks{The author M. Dai was partially supported by NSF grant DMS--1517583.}


\begin{abstract}
Considered here is the periodic initial-value probem for the 
 regularized long-wave (BBM) equation
\[u_t+u_x+uu_x-u_{xxt}=0.\]  
Adding to previous work in the literature, it is shown here that for any $s < 0$,
there is smooth initial data that is small in the $L_2$-based Sobolev
spaces $H^s$, but the solution emanating from it becomes arbitrarily large in 
arbitrarily small time.  This so called {\it norm inflation} result 
has as a consequence the previously determined conclusion that 
this problem is ill-posed in these negative-norm spaces.  

\end{abstract} 

\maketitle


\section{Introduction}\label{sec:intr}

This note derives from the paper \cite{BT} where it was shown that the initial-value problem\begin{equation}\label{bbm1}
\begin{split}
u_t+u&_x+uu_x-u_{xxt}=0,\\
&u(0,x)=u_0(x),
\end{split}
\end{equation}
 for the regularized  long-wave or BBM equation is globally well posed in the $L_2$--based Sobolev 
spaces $H^r(\R)$  provided $r \geq 0$.  In the same paper, it was shown that the  
map that takes initial data to solutions cannot be locally $C^2$ if $r < 0$.   This latter 
result suggests, but does not prove, that the problem \eqref{bbm1} is not well posed in $H^r$ for negative values of $r$.   Later, Panthee \cite{PA} showed that this solution map, 
were it to exist on all of $H^r(\R)$, could
not even be continuous, thus proving that the problem is 
ill posed in the $L_2$-based Sobolev spaces with negative index.  Indeed, he showed 
that there is a sequence of smooth initial data $\{\phi_n\}_{n = 1}^\infty$ such that 
$\phi_n \to 0$ in $H^r(\R)$ but the associated solutions, $\{u_n\}_{n = 1}^\infty$ have the property that $\|u(\cdot,t)\|_{H^r}$ is bounded 
away from zero for all small values of $t > 0$ and all $n \geq 1$.   

The BBM equation itself was initially put forward in \cite{P} and \cite{BBM} as an
approximate description of long-crested, surface water waves.   It is an 
alternative to the classical Korteweg-de Vries equation and has been shown to
be equivalent in that, for physically relevant initial data, the solutions of the two 
models differ by higher order terms on a long time scale
(see \cite{BPS1}.)   It  predicts the propagation of surface 
water waves pretty well in its range of validity \cite{BPS2}.  Finally, it is 
known rigorously to be a good approximation to solutions of the full, 
inviscid, water-wave problem by combining results in \cite{AABCW}, \cite{BCL} and \cite{L}
(see also \cite{L1}).  

It is our purpose here to show that in fact, for $r < 0$, the problem \eqref{bbm1} is not only 
not well posed, but features blow-up in the $H^r$--norm in arbitrarily 
short time.  This will be done in the context of the
periodic initial-value problem wherein $u_0$ is a periodic distribution 
lying in $H^r_{per}$ for some $r < 0$.   Similar results hold for $H^r(\R)$, but
are not explicated here.    

More precisely, it will be shown  that, for any given $r < 0$,  there is a sequence $\{u_0^n\}_{n = 1}^\infty$ of smooth initial
data such that $u_0^n \to 0$ in $H^r_{per}$ and  a  
sequence $\{T_n\}_{n=1}^\infty$ of positive times tending to  $0$ 
as $n \to \infty$ such that the corresponding solutions 
$\{u_n\}_{n=1}^\infty $ emanating from this initial data, whose existence is guaranteed by the periodic version \cite{C}
of the theory  for the initial-value problem,  are such that for $n = 1, 2, 3, \cdots$,
$$
\| u(\cdot,T_n) \|_{H^s_{per}} \geq n.
$$
This  insures in particular that the solution map $\mathcal{S}$ that associates solutions 
to initial data, which exists on $L_2$, cannot be extended continuously to all of 
$H^s_{per}$, thus reproducing Panthee's conclusion.  
Results of this sort go by the appellation {\it norm inflation} for obvious 
reasons.  The idea originated in the work of Bourgain and Pavlovi\'c \cite{BP} for the three-dimensional Navier-Stokes equation.  The method of construction there was applied to some other dissipative fluid equations by the second author and her collaborators, see \cite{DQS, CDnse, CDmhd}. It suggests that the method is generic as well as sophisticated. 
\vspace{.1cm} 

\centerline{\it Notation}  

The notation used throughout is standard.  For $r \in \R$, the 
 collection $\dot H^r_{per}$ is the homogeneous  space 
of $2\pi$-periodic distributions whose norm 
$$
\| f \|_r^2 \, = \, \sum_{k = 1}^\infty  k^{2r}(|f_k|^2 + |g_k|^2)
$$
 is finite.  Elements in   $\dot H^r_{per}$ all have mean zero over the period domain 
$[0,2\pi]$.  Here, the $\{f_k\}$ are the Fourier sine coefficients  and the 
  $\{g_k\}$ are the Fourier cosine coefficients of $f$.   Notice that 
$\dot H^0_{per}$  may be viewed simply as 
the $L_2$-functions on the period domain $[0,2\pi]$ with mean zero.
  If $X$ is any Banach 
space, the set $C([0,T];X)$ consists of the continuous functions from the 
real interval $[0,T]$ into $X$ with its usual norm.

\section{Norm inflation}

The principal result of our study is the following theorem.   

\begin{Theorem}\label{thm-NI}
Let $r<0$ by given.   Then there is  a sequence $\{ u^{j}_0\}_{j=1}^\infty$ 
of $\mathcal{C}^\infty$, periodic initial data such that 
$$
   u_0^{(j)} \to 0 \quad {\rm as} \quad j \to \infty
$$ 
in $\dot H^r_{per}$ and a sequence $\{T_j\}_{j=1}^\infty$ of positive times 
tending to zero as $j \to \infty$ such that  if $u_j(x,t)$ is the solution emanating 
from $u_0^{(j)}$, then 
$$
\| u(\cdot,T_j)\|_{\dot H^r_{per}} \, \geq \,  j
$$
for all $j = 1, 2, \cdots$.  
\end{Theorem} 
\pf  Fix  $s > 0$, let $r = -s$ and consider 
 a wavenumber $k_1 \in \mathbb{N}$ which, in due course, will be taken to be large.
  Let $k_2=k_1+1$, 
define $\bar{u}$ by 
$\bar{u}=\sin( {k_1x})+\sin ({k_2x})$ and consider the $2\pi$-periodic, men zero
 initial data $u_0= k_1^{\gamma}\bar u$ for \eqref{bbm1} where $\gamma > 0$ 
will be restricted presently.  
Of course, $u_0$ is smooth, so 
the theory developed in \cite{C} implies that a unique, global, smooth solution emanates from
this initial data.  Notice also that the solution preserves the property of having 
zero mean, so it lies in $C([0,T];\dot H^\rho_{per})$ for all $\rho \in \R$.   

Let $\varphi(D_x)$ be the Fourier multiplier operator given in terms of 
its Fourier transform  by $\widehat{\varphi(D_x)u}(\xi)=\frac{\xi}{1+\xi^2}\hat u(\xi)$. The equation (\ref{bbm1}) can be rewritten as 
\begin{equation} 
\begin{split}  \label{integral}
iu_t=&\varphi(D_x)u+\frac12 \varphi(D_x)\left(u^2\right),\\
&u(0,x)=u_0(x).
\end{split}
\end{equation}
Let $S(t)=e^{-it\varphi(D_x)}$ be the unitary group defining the evolution of 
the linear BBM equation. Then, Duhamel's principle allows 
the solution of (\ref{bbm1})-\eqref{integral} to be written in the form
\begin{equation}\label{decomposition}
u(x,t)=S(t)u_0(x)+u_1(s,t)+y(x,t)
\end{equation}
where
\begin{equation}\notag
u_1(x,t)=\frac12 \int_0^tS(t-\tau)\varphi(D_x)\big(S(\tau)u_0)\big)^2d\tau
\end{equation}
is the first order approximation of the nonlinear term in the 
differential-integral equation in
\eqref{integral}.  The function $y(x,t)$ is the remainder, which may 
be expressed implicitly in the sightly complicated, but useful form 
\begin{equation}\label{remainder}
y(x,t)=\int_0^tS(t-\tau)\varphi(D_x)\big[G_0(\tau)+G_1(\tau)+G_2(\tau)\big]\,d\tau
\end{equation}
with 
\begin{equation}\notag
\begin{split}
&G_0(\tau)=\frac12 u_1^2(\tau) + u_1(\tau)S(\tau)u_0,\\
&G_1(\tau)=u_1(\tau)y(\tau)+y(\tau)S(\tau)u_0,\\
&G_2(\tau)=\frac12 y^2(\tau),\\
\end{split}
\end{equation}
where the spatial dependence has been supressed for ease of reading.  
The strategy is to show that by choosing $k_1$ sufficiently large, $u_1$ becomes  large in a short time  in the space $\dot H^r_{per} = \dot H^{-s}_{per}$, while the error term $y$ remains under control in the same space.

In contrast to dissipative equations,  the linear 
 dispersion operator $S(t)$ only translates the wave, 
but does not change its magnitude; more precisely, for $k = 1, 2, \cdots$, 
\begin{equation} \label{dispersion}
S(t)\sin (kx)=\sin\left(kx-\frac{k}{1+k^2}t\right)\!, \;\; 
 S(t)\cos (kx)=\cos\left(kx-\frac{k}{1+k^2}t\right)\!.
\end{equation}
On the other hand, the operator $\varphi(D_x)$ both decreases the amplitude of 
its argument
 and adds rotation {\it viz.} 
\begin{equation}\label{phi}
\varphi(D_x)\sin kx=-i\frac{k}{1+k^2}\cos kx, \qquad \varphi(D_x)\cos kx=i\frac{k}{1+k^2}\sin kx.
\end{equation} 
It follows from this that $\varphi (D_x) $ vanishes on constant functions.  

It is clear that if $s > 0$, then
\begin{equation}\label{initial}
\begin{split}
&\|S(t)u_0\|_{-s}=\|u_0\|_{-s}\sim k_1^{\gamma-s}, \\  {\rm while} \quad \quad
&\|S(t)u_0\|_{0}=\|u_0\|_{0}\sim  k_1^{\gamma}.
\end{split}
\end{equation}
As we want the initial data to be small in $\dot H^{-s}_{per}$, $\gamma$ is restricted 
to the range $(0,s)$.
The formulas in \eqref{dispersion} imply 
\begin{equation}\notag
S(\tau) \bar u=\sin\left(k_1x-\frac{k_1}{1+k_1^2}\tau\right)+\sin\left(k_2x-\frac{k_2}{1+k_2^2}\tau\right)\!,
\end{equation}
so that 
\begin{equation}\notag
\begin{split}
\big[S(\tau) \bar u \big]^2=\, \frac12&\left[1-\cos\left(2k_1x-\frac{2k_1}{1+k_1^2}\tau\right)\right]+\frac12\left[1-\cos\left(2k_2x-\frac{2k_2}{1+k_2^2}\tau\right)\right]\\
&+\cos\left((k_1-k_2)x-\left(\frac{k_1}{1+k_1^2}-\frac{k_2}{1+k_2^2}\right)\tau\right)\\
&-\cos\left((k_1+k_2)x-\left(\frac{k_1}{1+k_1^2}+\frac{k_2
}{1+k_2^2}\right)\tau\right)\!.
\end{split}
\end{equation}
It then follows from \eqref{phi} that 
\begin{equation} \notag
\begin{split}
\frac12 \varphi(D_x)\big[S(\tau)\bar u\big]^2=&-\frac i4 \frac{2k_1}{1+4k_1^2}\sin\left(2k_1x-\frac{2k_1}{1+k_1^2}\tau\right)\\
&-\frac i4 \frac{2k_2}{1+4k_2^2}\sin\left(2k_2x-\frac{2k_2}{1+k_2^2}\tau\right)\\
&+\frac i2\frac{k_1-k_2}{1+(k_1-k_2)^2}\sin\left((k_1-k_2)x-\left(\frac{k_1}{1+k_1^2}-\frac{k_2}{1+k_2^2}\right)\tau\right)\\
&-\frac i2\frac{k_1+k_2}{1+(k_1+k_2)^2}\sin\left((k_1+k_2)x-\left(\frac{k_1}{1+k_1^2}+\frac{k_2}{1+k_2^2}\right)\tau\right)\\
\equiv & \,\, I_1+I_2+I_3+I_4.
\end{split}
\end{equation}

Consider now the function $\sin\left(kx - \omega t \right)$ and calculate as follows:
\begin{equation} \label{calculation}
\begin{split}
\int_0^t S(t-\tau)\sin(kx - \omega \tau)d\tau = \int_0^t \sin\Big(kx - \frac{k}{1+k^2}(t-\tau) - \omega \tau\Big) d\tau\\
= \left(\frac{k}{1+k^2} - \omega \right)^{-1} \left(\cos\Big(kx - \frac{k}{1+k^2} t \Big) 
   - \cos \Big(kx - \omega t \Big)\right)
\end{split}
\end{equation}
where use has been made of \eqref{dispersion}.

The latter
 formula, applied four times,  allows us to calculate $u_1$ explicitly, to wit, 
\begin{equation} \notag
\begin{split}
u_1=& k_1^{2\gamma}\int_0^tS(t-\tau)\big[I_1+I_2+I_3+I_4\big]\,d\tau\\
=&-\frac{ i k_1^{2\gamma}}{12}\frac{1+k_1^2}{k_1^2}\left[\cos\left(2k_1x-\frac{2k_1}{1+k_1^2}t\right)-\cos\left(2k_1x-\frac{2k_1}{1+4k_1^2}t\right)\right]\\
&-\frac{ik_1^{2\gamma}}{12}\frac{1+k_2^2}{k_2^2}\left[\cos\left(2k_2x-\frac{2k_2}{1+k_2^2}t\right)-\cos\left(2k_2x-\frac{2k_2}{1+4k_2^2}t\right)\right]\\
&+\frac{ik_1^{2\gamma}}{2}\frac{k_1-k_2}{1+(k_1-k_2)^2}\left[\frac{k_1}{1+k_1^2}-\frac{k_2}{1+k_2^2}-\frac{k_1-k_2}{1+(k_1-k_2)^2}\right]^{-1}\cdot \\
&\left[\cos\left((k_1-k_2)x-\left(\frac{k_1}{1+k_1^2}-\frac{k_2}{1+k_2^2}\right)t\right)-\cos\left((k_1-k_2)x-\frac{k_1-k_2}{1+(k_1-k_2)^2}t\right)\right]\\
&-\frac{ik_1^{2\gamma}}{2}\frac{k_1+k_2}{1+(k_1+k_2)^2}\left[\frac{k_1}{1+k_1^2}+\frac{k_2}{1+k_2^2}-\frac{k_1+k_2}{1+(k_1+k_2)^2}\right]^{-1}\cdot \\
&\left[\cos\left((k_1+k_2)x-\left(\frac{k_1}{1+k_1^2}+\frac{k_2}{1+k_2^2}\right)t\right)-\cos\left((k_1+k_2)x-\frac{k_1+k_2}{1+(k_1+k_2)^2}t\right)\right].
\end{split}
\end{equation}
A study of the various constants appearing above reveals that, up to absolute constants, 
 \begin{equation}\notag
 \begin{split}
 u_1\sim&-ik_1^{2\gamma}\left[\cos\left(2k_1x-\frac{2k_1}{1+k_1^2}t\right)-\cos\left(2k_1x-\frac{2k_1}{1+4k_1^2}t\right)\right]\\
 &-ik_1^{2\gamma}\left[\cos\left(2k_2x-\frac{2k_2}{1+k_2^2}t\right)-\cos\left(2k_2x-\frac{2k_2}{1+4k_2^2}t\right)\right]\\
 &+i k_1^{2\gamma}\left[\cos\left(x-\left(\frac{k_1}{1+k_1^2}-\frac{k_2}{1+k_2^2}\right)t\right)-\cos\left(x-\frac {t}{2}\right)\right]\\
 &-i k_1^{2\gamma}\Bigg[\cos\left((k_1+k_2)x-\left(\frac{k_1}{1+k_1^2}  +\frac{k_2}{1+k_2^2}\right)t\right)  \\  
& \hspace{.95cm} -\cos\left((k_1+k_2)x-\frac{k_1+k_2}{1+(k_1+k_2)^2}t\right)\Bigg].
 \end{split}
 \end{equation}
as $k_1$ becomes large.  Since 
$$
\big|\cos(kx - \omega_1 t)  -\cos(kx - \omega_2 t)\big| \, \leq \,  |\omega_1 - \omega_2|t, 
$$
straightforward calculations show that the first, second and fourth terms above 
are uniformly small compared to the third term, for large values of $k_1$.  Indeed, 
they are all of order $ k_1^{2\gamma-1}t$, whereas the third term is of order 
 $k_1^{2\gamma} t$.

It follows from this that for all $t \geq 0$, 
\begin{equation} \label{u1}
\begin{split}
\|u_1(t,\cdot)\|_{-s}\sim &\, k_1^{2\gamma} t \quad {\rm and \; likewise}\\
\|u_1(t,\cdot)\|_0\sim& \,k_1^{2\gamma} t.
\end{split}
\end{equation}
Thus, by taking $k_1$ large, the $\dot H^{-s}_{per}$-norm of $u_1$ 
can be made  as big as we like.

  As mentioned earlier, an estimate of the error term $y$ is needed to complete the argument.   It will in fact be shown that $y$ 
is even bounded  in $L_2$, let along $\dot H^{-s}_{per}$. 
 To this end, use is made of one 
of a periodic version of one the  bilinear estimates in \cite{BT}.
\begin{Lemma} \label{le:bil}
Let $u,v\in H^q_{per}$ with $q\geq 0$. Then 
\begin{equation}\label{bilinear}
\|\varphi(D_x)(uv)\|_{q}\lesssim \|u\|_q\|v\|_{q}
\end{equation}
where the implied constant only depends upon $q$.   
\end{Lemma}

The proof of this  result is the same as the proof of Lemma 1 in \cite{BT}, with sums replacing integrals.   

\vspace{.1cm}

Introduce  the abbreviation $X_T$ for $C([0,T];L^2)$
for ease of reading.  The value of $T > 0$ will be specified momentarily.  
It follows from \eqref{bilinear} and the implicit relationship \eqref{remainder} 
for the remainder $y$ that
\begin{equation}\label{final}
\begin{split}
\|y\|_{X_T}\lesssim &T \|u_1\|_{X_T}^2+T\|S(t)u_0\|_{X_T}\|u_1\|_{X_T}+T\|u_1\|_{X_T}\|y\|_{X_T}\\
&+T\|S(t)u_0\|_{X_T}\|y\|_{X_T}+T\|y\|_{X_T}^2\\
\lesssim &T^3k_1^{4\gamma}+T^2k_1^{3\gamma}+\Big(k_1^{2\gamma} T^2
+ k_1^{\gamma}T \big)\|y\|_{X_T} +T\|y\|_{X_T}^2 \\
=& \mathcal{A} + \mathcal{B} \mathcal{Y} + T \mathcal{Y}^2,
\end{split}
\end{equation}
where $\mathcal{Y} = \mathcal{Y}(T) = \|y\|_{X_T}$.  As $y \in C([0,M];L_2)$ for 
all $M > 0$, it follows that $\mathcal{Y}(T)$ is a continuous function of $T$.  
Moreover, $\mathcal{Y}(0) = 0$. 

Choose $T_0 = k_1^{-\mu\gamma}$, where $\mu > \frac32$. With this choice, we see that for $T \leq T_0$,  
$$
\mathcal{A} =O\big( k_1^{\gamma(4-3\mu)} + k_1^{\gamma(3-2\mu)}\big) \quad {\rm and} \quad \mathcal{B} = O\big(k_1^{2\gamma(1-\mu)} + k_1^{\gamma(1-\mu)} \big),
$$  
as $k_1 \to \infty$ and all the exponents are negative.  

Choose $k_1$ large enough that 
$\mathcal{B} < \frac12$ and $T$ and $\mathcal{A}$ are both small.   It follows in 
this circumstance that
the quadratic polynomial
$$
p(z) = \mathcal{A} + (\mathcal{B} - 1)z + Tz^2
$$ 
has two  positive roots, the smaller of which is denoted $\underbar z$ and the
larger $\bar z$.  Of course, $p(z) < 0$ for $z \in (\underbar z, \bar z)$.  

The inquality \eqref{final} may be expressed as 
$$
p\big(\mathcal{Y}(T)\big) \, \geq \, 0.
$$
As $\mathcal{Y}(T)$ is continuous and $\mathcal{Y}(0) = 0$, it follows that
$\mathcal{Y}(T) \leq \underbar z$ for all $T \in [0,T_0]$.   
For $k_1$ large, $T_0 < 1$. When combined with the fact that $\mathcal{B} < \frac12$, 
it is readily deduced that
$$
            \underbar z \, \leq \, 4\mathcal{A}, \quad {\rm whence} \quad \mathcal{Y}(T) \, \leq \, 4\mathcal{A},
$$
thus assuring that the remainder $y(\cdot,t)$ is indeed uniformly 
bounded in $\dot H^{-s}_{per}$ 
for $t \leq T_0$ and large choices of $k_1$.     

Taking a suitably chosen, increasing sequence $\{k_1^{(j)}\}_{j=1}^\infty$ 
of wavenumbers 
for which 
$$
\lim_{j \to \infty} k_1^{(j)} = +\infty,
$$ 
and with the indicated choices of $\gamma$ and $\mu$, \eqref{initial} assures the initial data tends to zero in $\dot H^{-s}_{per}$.  The 
decomposition \eqref{decomposition} together with \eqref{initial}, \eqref{u1} 
and the bound just obtained on $y$ then implies that the solutions $u_j$ 
blow up at times $T_j = (k_1^{(j)})^{-\mu\gamma}$. The latter tend to zero as 
$j \to \infty$ since $\mu$ and $\gamma$ are both positive.
   This completes the proof of the theorem.


{}


\begin{thebibliography}{References}

\bibitem{AABCW}  
\newblock A.A. Alazman, J.P. Albert, J.L. Bona, 
M. Chen and J. Wu.  
\newblock  {\em Comparisons between the BBM equation and a
Boussinesq system}.
\newblock Advances Differential. Eq. {\bf 11}  
\newblock  (2006) 121--166.

\bibitem{ABN}
\newblock D. Ambrose, J.L. Bona, and D. Nicholls.
\newblock {\em On Ill-posedness of truncated series models for water waves}.
\newblock   Proc.  Royal Soc. London, Series A {\bf 470} 
\newblock  (2014) 1--16.

\bibitem{BBM} T.B. Benjamin,  J.L. Bona and J.J. Mahony.  
\newblock {\em Model equations for long waves in nonlinear dispersive media}, 
\newblock    Philos. Trans. Royal Soc. London Series A {\bf 272}
\newblock (1972) 47--78.

\bibitem{BCL}  J.L. Bona, T. Colin and D. Lannes.
\newblock  {\em Long wave approximations for water waves},
\newblock  Archive Rat. Mech. Anal. {\bf 178}
\newblock  (2005)  373--410.

\bibitem{BPS2}  J.L. Bona, W.G. Pritchard and L.R. Scott.  
\newblock {\em An evaluation of a model equation for water waves},
 \newblock Philos. Trans. Royal Soc. London Series A {\bf 302}
\newblock (1981) 457--510.

\bibitem{BPS1}  J.L. Bona, W.G. Pritchard and L.R. Scott.  
\newblock {\em A comparison of solutions of two model equations for long waves},
\newblock In Lectures in Applied Mathematics {\bf 20} 
\newblock  (ed. N. Lebovitz)   American Mathematical Society: Providence 
\newblock (1983) 235--267. 


\bibitem{BT}   J.L. Bona and N. Tzvetkov.  
\newblock {\em Sharp well-posedness results for  the BBM equation}, 
\newblock   {Discrete \&
Continuous Dynamical Systems, Series A {\bf 23} }
\newblock (2009) 1241--1252. 

\bibitem{BP}
\newblock .
\newblock {\em Ill-posedness of the Navier-Stokes equations in a critical space in 3D}.
\newblock   Journal of Functional Analysis,  {\bf 255} 
\newblock  (2008) 2233--2247.

\bibitem{C}  H. Chen.  
\newblock  {\em Periodic initial-value problem for the BBM-equation}, 
\newblock Computers and Mathematics with Applications, Special Issue on “Computational Methods in Analysis” {\bf 48} (2004) 1305--1318.

\bibitem{CDmhd}
\newblock A. Cheskidov and M. Dai.
\newblock {\em Norm inflation for generalized Magneto-hydrodynamic system}.
\newblock  {Nonlinearity, {\bf 28}}
\newblock (2015)  129--142.

\bibitem{CDnse}
\newblock A. Cheskidov and M. Dai.
\newblock {\em Norm inflation for generalized Navier-Stokes equations}.
\newblock   {Indiana University Mathematics Journal, {\bf 63}}
\newblock (2014), No. 3 : 869--884.

\bibitem{DQS}
\newblock M. Dai, J. Qing, and M. Schonbek.
\newblock {\em Norm inflation for incompressible Magneto-hydrodynamic system in $\dot{B}_\infty^{-1,\infty}$}.
\newblock   {Advances in Differential Equations, {\bf 16}}
\newblock (2011), No. 7-8, 725--746.

\bibitem{L}  D. Lannes.
\newblock {\em Well-posedness of the water-waves equations},
\newblock  J. American Math. Soc. {\bf 18}
\newblock  (2005)  605--654.

\bibitem{L1}  D. Lannes.
\newblock  The water waves problem: mathematical analysis and asymptotics,
\newblock  {\em Mathematical Surveys and Monographs {\bf 188} }
\newblock  American Math. Soc.: Providence (2013).

\bibitem{PA}   M. Panthee, 
\newblock {\em On the ill-posedness result for the BBM equation} 
\newblock   {Discrete \&
Continuous Dynamical Systems {\bf 30} }
\newblock (2011) 253--259. 

\bibitem{P}  D.H. Peregrine, 
\newblock {\em Calculations of the development of an undular bore}, 
\newblock   {J. Fluid Mech. {\bf 25}}
\newblock (1966) 321--330 . 

\end{thebibliography}
\end{document}